\documentclass[9pt]{amsart}

\usepackage{color}
\usepackage{graphicx}

\newtheorem{theorem}{Theorem}[section]

\newtheorem{corollary}[theorem]{Corollary}
\newtheorem{algorithm}[theorem]{Algorithm}
\newtheorem{claim}[theorem]{Claim}

\theoremstyle{definition}
\newtheorem{definition}[theorem]{Definition}

\theoremstyle{remark}

\numberwithin{equation}{section}

\newcommand{\numsixclusters}{435696}
\newcommand{\numsevenclusters}{1154}
\newcommand{\highlight}[1]{#1} %% {\textcolor{red}{#1}}
\newcommand{\highlightt}[1]{#1} %% {\textcolor{green}{#1}}

\def\Ps{\mathcal{P}}
\def\di{m}

\begin{document}

\title{Constructing $7$-clusters}

\author{Sascha Kurz}
\address{Sascha Kurz, Department of Mathematics, Physics and Informatics, University of Bayreuth, Bayreuth, Germany}
\email{sascha.kurz@uni-bayreuth.de}

\author{Landon Curt Noll}
\address{Landon Curt Noll, Cisco Systems, San Jose California, USA}
\email{ncluster-mail@asthe.com}

\author{Randall Rathbun}
\address{\highlightt{Randall Rathbun, Green Energy Technologies, LLC, Manning, Oregon, USA}}
\email{randallrathbun@gmail.com}

\author{Chuck Simmons}
\address{Chuck Simmons, Google, Mountain View, California, USA}
\email{csimmons@google.com}



\begin{abstract}
  A set of $n$-lattice points in the plane, no three on a line and no four
  on a circle, such that all pairwise distances and all coordinates
  are integral is called an $n$-cluster (in $\mathbb{R}^2$). We determine the
  smallest existent $7$-cluster with respect to its diameter. Additionally
  we provide a toolbox of algorithms which allowed us to computationally
  locate over 1000 different $7$-clusters, some of them having 
  huge integer edge lengths. On the way, we exhaustively determined
  all Heronian triangles with largest edge length up to $6\cdot 10^6$.
\end{abstract}

\maketitle

\section{Introduction}
Point sets with pairwise rational or integral distances have been studied for a
long time; see e.g.\ \cite{48.0667.04,kummer}. For brevity we will call those point sets rational or integral.
Nevertheless only a few theoretical results are known and integral point sets seem to be unexpectedly
difficult to construct. On the other hand there is the famous open problem, asking for a dense set in
the plane such that all pairwise Euclidean distances are rational, posed by Ulam in 1945, see
e.g.\ \cite{0086.24101}. Till now we only know that one can easily construct rational points sets
which are either dense on a line or a circle, see e.g., \cite[Sec.~5.11]{1086.52001} or \cite{dense_subset}.
In \cite{1209.52009} the authors have shown that no irreducible algebraic curve other than a line or a circle
contains an infinite rational set. Thus if Ulam's question admits a positive answer the corresponding point
set has to be very special.

\highlight{Almering \cite{almering1963rational} established, that for a given triangle with rational side lengths,
the set of points with rational distances to the three vertices, is dense in the plane of the triangle. Berry \cite{berry1992points}
relaxed the conditions to one rational side length and the other two side lengths being a square root of a rational number.
More general considerations can be found in the recent preprint \cite{barbara2013rational}. So far no such result
is known for a quadrilateral with pairwise rational distances. With just one distance missing 
\highlightt{Dubickas states in \cite{dubickas2012some}}
that every $n\ge 3$ points in $\mathbb{R}^2$ can be slightly perturbed to a set of $n$ points in $\mathbb{Q}^2$ such
that at least $3(n-2)$ of the mutual distances are rational. Declaring which of the mutual distances has to be rational
can be modeled as a graph. Classes of \textit{admissible} graphs have been studied e.g.\ in \cite{benediktovich2013rational,
geelen2008straight}.}

Given a finite rational point set, we can of course convert it into an integral
point set by rescaling its edge lengths with the least common multiple of their respective denominators.\footnote{As
shown in \cite{ErdoesAnning1,ErdoesAnning2} each infinite integral point set is located on a line.} Thus for each
finite number $n$ one can easily construct an integral point set consisting of $n$ points where (almost) all points
are located on a circle. Constructions of finite integral point sets where  $n-1$ or $n-2$ points are located on a
line are e.g.\ given in \cite{1265.52018}. To this end several authors, including Paul Erd\H{o}s \cite[Problem D20]{UPIN},
ask for integral points sets in general position, meaning that no three points are on a line and no four points are on a
circle. These objects seem to be rather rare or at the very least hard to find. For $n=6$ points a few general constructions
for integral points sets in general position are known \cite{hab_kemnitz}. The only two published examples of $7$-point
integral point sets in general position are given in \cite{1145.52010}. Independently and even earlier, in May 2006
Chuck Simmons and Landon Curt Noll found\footnote{cf. http://www.isthe.com/chongo/tech/math/n-cluster/} even more
restricted configurations. At that time the \textit{smallest} one has (integral) coordinates
\begin{eqnarray*}
  \!\!\!\!\!\!\! &&(0,0) (327990000,0) (238776720,118951040) (222246024,-103907232)\\
  \!\!\!\!\!\!\! &&(243360000,21896875) (198368352,50379264) (176610000,-94192000)
\end{eqnarray*}

Aiming at $n$-point integral point sets in general position, for especially $n=6$, Noll and Bell \cite{cluster} additionally
required that also the coordinates have to be integral and called those structures $n_2$-clusters, or when the restriction
to the dimension\footnote{The notion of an integral point set can be easily generalized to arbitrary dimensions $m$. The
term \textit{general position} then has the meaning that no $m+1$ points are contained in a hyperplane and no $m+2$ points
are contained in a hypersphere, \highlightt{see e.g.\ \cite{cluster}.}} is clear from the context, $n$-clusters. Using a 
computer search the authors found
91~non-similar $6$-clusters, where the respective greatest common divisor of their corresponding edge lengths is one,
but no $7$-clusters.\footnote{Independently also Randall Rathbun found the first few $6$-clusters.} Using a slightly
improved version and lots of computing time Simmons and Noll in  2006 found the first $7$-clusters and extended there list
to twenty-five $7$-clusters in 2010.

The aim of this paper is to present a set of sophisticated algorithms in order to construct $n$-clusters for $n\ge 7$.
Using an exhaustive search we were able to determine, with respect to its diameter, the smallest $7$-cluster and provide
heuristic methods to produce more than 1000 non-similar $7$-clusters. Unfortunately so far no $8$-cluster turned up.
So the hunt for an integral octagon in general position or even an $8$-cluster is still open. In this context we mention 
the Erd\H{o}s/Noll infinite-or-bust $n_m$-cluster conjecture: For any dimension $m > 1$,  and any number of points $n > 2$,
there exists either $0$ or an infinite number of primitive $n_m$-clusters.

In Section~\ref{sec_basic} we summarize the known theory on integral point sets and in Section~\ref{sec_heronian_triangles}
we go into the algorithmic details how to generate large lists of Heronian triangles. Section~\ref{sec_exhaustive_n_cluster}
is devoted to exhaustive searches for $n$-clusters up to a given diameter. Here the idea is to combine $n$-clusters that share
a common $n-1$-cluster. Allowing the containment of similar $n-1$-clusters, i.e.\ a scaled version, is the idea behind
Section~\ref{sec_combine}. Our most successful algorithmic approach is presented in Section~\ref{sec_triangle_extension}. Since
the basic operations of our algorithms have to be performed quite often, we present low level details in Section~\ref{sec_low_level_details}.
A theoretically interesting algorithm, based on circle inversion, is presented in Section~\ref{sec_circle_inversion}. 
\highlight{Methods to extend a given triangle with rational side length by a forth point are studied in
Section~\ref{sec_fourth_point}.} Since
almost all of our presented algorithms depend on a selection of Heronian triangles, which may not be too large due to
computational limits, we present ways to select Heronian triangles from larger sets in Section~\ref{sec_promising}. Our
computational observations are summarized in Section~\ref{sec_observations}. We present our computational results in
Section~\ref{sec_results} before we draw a conclusion in Section~\ref{sec_conclusion}.

\section{Basic results and notation}
\label{sec_basic}

\begin{definition}
  An integral point set $\mathcal{P}$ is a set of points in the plane that are not all located on a line such that
  the pairwise differences are all integers.
\end{definition}

We remark that integral point sets can easily be defined in arbitrary dimensions, see e.g.\ \cite{1088.52011,1135.51305} but
the present paper is restricted to the two-dimensional case.

One of the first question arising when dealing with integral point sets is how to represent them. Of course one may list
a coordinate representation. One example of such a representation is in the introduction. Another way is to provide a
table of the pairwise distances from which a coordinate representation can easily be computed. For the example from
the introduction we have the following distance table:
$$
{\tiny
  \left(\begin{array}{cccccccc}
    0 & 327990000 & 266765200 & 245336520 & 244343125 & 204665760 & 200158000 \\
    327990000 & 0 & 148688800 & 148251480 & 87416875 & 139067760 & 178292000 \\
    266765200 & 148688800 & 0 & 223470520 & 97162325 & 79592240 & 222024000 \\
    245336520 & 148251480 & 223470520 & 0 & 127563605 & 156123240 & 46658680 \\
    244343125 & 87416875 & 97162325 & 127563605 & 0 & 53249365 & 133911125 \\
    204665760 & 139067760 & 79592240 & 156123240 & 53249365 & 0 & 146199440 \\
    200158000 & 178292000 & 222024000 & 46658680 & 133911125 & 146199440 & 0 \\         
  \end{array}\right)
}  
$$
Given a matrix of distances one can decide whether there exists a set of vertices
in the $m$-dimensional Euclidean space $\mathbb{R}^m$ attaining those distances
based on a set of inequalities and equations involving the so-called Cayley-Menger determinants
\cite{paper_characteristic,menger}. 

\begin{definition}
  If $\Ps$ is a point set in $\mathbb{R}^\di$ with vertices $v_0,v_1,\dots,v_{n-1}$ and 
  $C=(d_{i,j}^2)$ denotes the $n\times n$ matrix given by $d_{i,j}^2=\Vert v_i-v_j\Vert_2^2$ 
  the Cayley-Menger matrix $\hat{C}$ is obtained from $C$ by bordering $C$ with a top 
  row $(0,1,1,\dots,1)$ and a left column $(0,1,1,\dots,1)^T$.
\end{definition}

\begin{theorem}{\textbf{(Menger \cite{menger})}}
  \label{theorem_menger}
  A set of vertices $\{v_0,v_1,\dots,v_{n-1}\}$ with pairwise distances $d_{i,j}$ is realizable 
  in the Euclidean space $\mathbb{R}^{\di}$ if and only if for all subsets $\{i_0,i_1,\dots,i_{r-1}\}
  \subset\{0,1,\dots,n-1\}$ of cardinality $r\le \di+1$,
  $$
    (-1)^{r}CMD(\{v_{i_0},v_{i_1},\dots,v_{i_{r-1}}\})\ge 0,
  $$
  and for all subsets of cardinality $\di+2\le r\le n$,
  $$
    (-1)^{r}CMD(\{v_{i_0},v_{i_1},\dots,v_{i_{r-1}}\})=0\,.
  $$
\end{theorem}

Thus it is possible to deal with integral point sets by storing their pairwise distances only but often
it is computationally cheaper to use coordinate representations which are easy to compute. As remarked in 
the introduction we are interested in integral point sets in the Euclidean plane $\mathbb{R}^2$ with some
additional properties.

\begin{definition}
  An (plane) integral point set is in general position if no three points are on a line and no four points
  are on a circle.
\end{definition}

The condition on the arrangement of points can easily be generalized to higher dimension and also be expressed
using the Cayley-Menger determinants, see e.g.\ \cite{1088.52011,paper_characteristic}. For the plane it suffices
to check the triangle inequality in order to discover three collinear points. Checking the condition of Ptolemy's
theorem, one can easily discover four points on a circle.

\begin{definition}
  An $n$-cluster is a plane integral point set in general position that consists of $n$ points such
  that there exists a representation using integer coordinates, \highlightt{i.e., lattice points}.
\end{definition}

Fortunately we do not have to deal with the constraint of integral coordinates, but we have to go far afield:
The area $A_\Delta(a,b,c)$ of a triangle with side lengths $a$, $b$, $c$ is given by 
$$
  A_\Delta(a,b,c)=\frac{\sqrt{(a+b+c)(a+b-c)(a-b+c)(-a+b+c)}}{4}
$$
due to the Heron formula. If the area is non-zero, we can uniquely write  $A_\Delta(a,b,c)=q\sqrt{k}$ with a
rational number $q$ and a square-free integer $k$. The number $k$ is called the characteristic $\Delta$ of the
triangle with side lengths $a$, $b$, $c$. Kemnitz \cite{hab_kemnitz} has shown that each non-degenerate triangle
of an integral point set has the same characteristic, which was also generalized to arbitrary dimensions in
\cite{paper_characteristic}. Since triangles with integral coordinates have a rational area, \highlightt{see e.g.\ 
Pick's theorem,} the triangles of an $n$-cluster all have to have a characteristic of $1$.

We now argue that the opposite is also true. Given an integer sided triangle with characteristic $1$ we can
easily determine a representation using rational coordinates, see e.g.\ \cite{paper_characteristic}. Due to
Fricke \cite{Fricke}, see also \cite{lattice_tetrahedra,0989.11014}, each integral point set in the plane which
has a representation in rational coordinates has a representation in integral coordinates. Thus there is no need
to explicitly search for integral coordinates for $n$-clusters. One just needs to check that all pairwise distances
are integral and that at least one contained non-degenerate triangle has characteristic $1$ or, equivalently, 
that it has a representation in rational coordinates.

A Heronian triangle is a triangle with integer side lengths and area\footnote{Some authors allow the side lengths and
the area to be rational and remark that all quantities can be easily rescaled to be integers.}. Due to the formula
for $A_\Delta(a,b,c)$ for an integer sided triangle with characteristic $1$, the area is rational and may in principle be
non-integral. Nevertheless one may consider the cases of the side lengths modulo $8$  (see \cite{Chisholm2006153}) and
conclude that such triangles have to be integral. We summarize these findings in:
\begin{corollary}
  Given a non-degenerate triangle $T$ with integer side lengths then the following statements are equivalent:
  \begin{itemize}
    \item[(a)] $T$ has characteristic $1$
    \item[(b)] $T$ has rational area
    \item[(c)] $T$ has integral area, i.e.\ $T$ is Heronian
  \end{itemize}
\end{corollary}

Thus Heronian triangles are the basic building blocks of $n$-clusters and we will consider algorithms how to generate
them in the next section.

\bigskip

In the introduction we have spoken of \textit{the smallest} cluster. So in order to have a measure of the \textit{size}
of an $n$-cluster or more generally an integral point sets we denote the largest distance between two points
as its diameter.  If we perform \highlightt{an} exhaustive search in the following we will always have to impose a limit on the
maximum diameter. We remark that other metrics are possible too, but most of them can be bounded by constants
in terms of the maximum diameter.

\bigskip

Given an $n$-cluster we can obviously construct an infinite sequence of non-isomorphic $n$-clusters by rescaling
the clusters by integers $2,3,\dots$. We call those $n$-clusters similar and are generally interested in lists
of non-similar $n$-clusters. To this end we call a given $n$-cluster primitive if its edge lengths do not have
a common factor larger then $1$. As argued before dividing the edge lengths of a given integral point set
by the greatest common divisor does not destroy the property of admitting integral coordinates.

Applying this insight to the example given in the introduction we observe that the greatest common divisor
of the edge length is $145$. Thus dividing all edge lengths gives the following distance matrix:
$$
  \left(
  \begin{array}{rrrrrrr}
    0 & 2262000 & 1839760 & 1691976 & 1685125 & 1411488 & 1380400 \\
    2262000 & 0 & 1025440 & 1022424 & 602875 & 959088 & 1229600 \\
    1839760 & 1025440 & 0 & 1541176 & 670085 & 548912 & 1531200 \\
    1691976 & 1022424 & 1541176 & 0 & 879749 & 1076712 & 321784 \\
    1685125 & 602875 & 670085 & 879749 & 0 & 367237 & 923525 \\
    1411488 & 959088 & 548912 & 1076712 & 367237 & 0 & 1008272 \\
    1380400 & 1229600 & 1531200 & 321784 & 923525 & 1008272 & 0 \\
  \end{array}
  \right)
$$
This $7$-cluster has a diameter of $2\,262\,000$, which is smallest possible as verified in
Section~\ref{sec_exhaustive_n_cluster}. A coordinate representation is given by
\begin{eqnarray*}
  \!\!\!\!\!\!\!\!\!&&(0,0)(374400,-2230800)(1081600,-1488240)(-453024,-1630200)\\
  \!\!\!\!\!\!\!\!\!&&(426725,-1630200)(569088,-1291680)(-439040,-1308720) 
\end{eqnarray*}

\section{Generation of Heronian triangles}
\label{sec_heronian_triangles}
The conceptually simplest algorithm to exhaustively generate all Heronian triangles up to a given diameter
is to loop over all non-isomorphic integer triangles and to check whether the area is integral:
\begin{algorithm}
  {\textbf{(Exhaustive generation of Heronian triangles)}\\}
   for $a$ from $1$ to $n$\\ 
  \hspace*{3mm}for $b$ from $\lceil\frac{a+1}{2}\rceil$ to $a$\\
  \hspace*{6mm}for $c$ from $a+1-b$ to $b$\\ 
  \hspace*{9mm}if $\sqrt{(a+b+c)(a+b-c)(a-b+c)(-a+b+c)}\in\mathbb{N}$\\
  \hspace*{9mm}then output $a$, $b$, and $c$ 
\end{algorithm}

Assuming that the check in the last-but-one line can be performed in constant time, this algorithm
has time complexity $\Theta(n^3)$. Two $O(n^{2+\varepsilon})$ algorithms, where $\varepsilon>0$ is
arbitrary, have been given in \cite{herontriangles}.

Complete parameterizations have been known for a long time, i.e.\ the Indian mathematician Brahmagupta 
(598-668 A.D.) who gives, see e.g.\ \cite{45.0283.11,herontriangles}, the parametric solution 
\begin{eqnarray}
  a&=&\frac{p}{q}h(i^2+j^2)\nonumber\\
  b&=&\frac{p}{q}i(h^2+j^2)\nonumber\\
  c&=&\frac{p}{q}(i+h)(ih-j^2)\nonumber
\end{eqnarray} 
for positive integers $p$, $q$, $h$, $i$, and $j$ fulfilling $ih>j^2$ and $gcd(p,q)=gcd(h,i,j)=1$.

Due to the presence of the denominators $q$ this parameterization is not well compatible with
restrictions on the maximum diameter. On the other hand we can easily generate primitive, meaning
that the side lengths have no common factor, Heronian triangles by looping over all feasible 
triples $(h,i,j)$ below a suitable upper bound, setting $p$ to $1$ and choosing $q$ case dependent
such that $gcd(a,b,c)=1$. Using this approach we can quickly generate a huge amount of primitive 
Heronian triangles, but on the other hand may get those with small diameters rather late, compared
to the upper bound on $h,i,j$, and have to face the fact that the same primitive Heronian triangle
may be generated several times.

For the purpose of this paper we use a different exhaustive algorithm to generate all primitive Heronian
triangles up to a prescribed diameter. Given a triangle with side lengths $a$, $b$, and $c$ we have
$\cos\alpha=\frac{b^2+c^2-a^2}{2bc}$ and $\sin\alpha=\frac{2A_\Delta(a,b,c)}{bc}$. For a Heronian
triangle $\sin\alpha$ and $\cos\alpha$ are rational numbers so that also 
$\tan \frac{\alpha}{2}=\frac{\sin\alpha}{1+\cos\alpha}\in\mathbb{Q}$.
Thus there are coprime integers $m,n$ satisfying $\tan \frac{\alpha}{2}=\frac{n}{m}$. With these
parameters we obtain
$$
  \cos\alpha=\frac{1-\tan^2\frac{a}{2}}{1+\tan^2\frac{a}{2}}=\frac{m^2-n^2}{m^2+n^2},
$$
where $\gcd(m^2-n^2,m^2+n^2)\in\{1,2\}$. From the other representation of $\cos \alpha$ we can then
conclude that $m^2+n^2$ divides $2bc$. So given two integral side lengths $b$ and $c$ of a Heronian triangles,
we can determine all possibilities for $m^2+n^2$, then determine all possibilities for $m$ and $n$, and
finally determine all possibilities for the third side $a$:

\begin{algorithm}
  {\textbf{(Find the third side)}\\}
   loop over all divisors $k$ of $2bc$\\ 
  \hspace*{3mm}loop over all solutions $(m,n)$ of $m^2+n^2=k$\\
  \hspace*{6mm}solve $\frac{b^2+c^2-a^2}{2bc}=\frac{m^2-n^2}{m^2+n^2}$ for $a$\\ 
  \hspace*{9mm}if $a\in \mathbb{Q}$ and the triangle inequalities are satisfied for $(a,b,c)$\\
  \hspace*{9mm}then output $a$ 
\end{algorithm} 

So in order to determine all primitive Heronian triangles up to diameter $N$ we have to loop over all
coprime pairs $(b,c)$ with $N\ge b\ge c\ge 1$ and apply the above algorithm to determine $a$. Given $a$
we can check whether $a,b,c$ are coprime, $a\le N$, and $a\ge b$, $a\in\mathbb{N}$ (to avoid isomorphic duplicates). 

In this context the maximum diameter $n$ has to be limited to a few millions so that we can easily determine
the prime factorizations of all integers at most $n$ in a precomputation. Given those data we can 
quickly determine the prime factorization of $2bc$ and loop over all divisors without any additional testing. 

Next we want to describe the set of solutions of $m^2+n^2=k$ and assume that 
$$
  k=2^h\cdot q_1^{i_1}\dots q_s^{i_s}\cdot p_1^{j_1}\dots p_t^{j_t},
$$
where the $q_l$ are primes congruent to $3$ modulo $4$ and the $p_l$ are primes congruent to $1$ modulo $4$.
If any of the $i_l$ is odd, then no integer solution to $m^2+n^2=k$ exists. Otherwise each solution can be written
as $(m,n)=\lambda\cdot(\tilde{m},\tilde{n})$, where $\lambda=2^{\lfloor h/2\rfloor}\cdot q_1^{i_1/2}\dots q_s^{i_s/2}$
and $\tilde{m}^2+\tilde{n}^2=k/\lambda^2=:\tilde{k}$, i.e.\
$$
  \tilde{k}=2^{\tilde{h}}\cdot p_1^{j_1}\dots p_t^{j_t},
$$
where $\tilde{h}\le 1$. Due to the formula $(x_1^2+x_2^2)(y_1^2+y_2^2)=(x_1y_1+x_2y_2)^2+(x_1y_2-x_2y_1)^2$
and the unique factorization of the Gaussian integers $\mathbb{Z}[i]$ it suffices to combine the solutions the
problem, where $\tilde{k}$ is a prime power. Ignoring signs for $\tilde{k}=2$ the unique solution is given by $1^2+1^2=2$.
Ignoring signs and order then there is a unique solution for $u^2+v^2=p$ once $p$ is equivalent to $1$ modulo $4$.
Again ignoring signs and order, for prime powers the set of solutions $x^2+y^2=p^j$ is given by
$x+yi=(u+vi)^{l}(u-vi)^{j-l}$, where $0\le l\le j/2$. Thus it remains to determine a solution of $u^2+v^2=p$, which can be
done by the Hermite-Serret algorithm, which first determines an integer $z$ satisfying $z^2\equiv i\pmod p$, using that 
$w^{\frac{p-1}{2}}\equiv -1\pmod p$ for each quadratic nonresidue $w$, and then applies the Euclidean algorithm on $(p,w)$
to determine $(u,v)$. See \cite{hermite,serret1848theoreme} for the original sources and \cite{brillhart1972note} for an
improved algorithm. The just sketched algorithm for the generation of all Heronian triangles up to
diameter $n$ runs in $O(n^{2+\varepsilon})$ time, where $\varepsilon>0$ is arbitrary.

Using this algorithm we have exhaustively generated all primitive Heronian triangles up to diameter $6\cdot 10^6$.
They are available for download at \cite{HP}. Having the data at hand we
have computed an approximate counting function which fits best for a given type of functions. Let $count(x)$ denote
the number of primitive Heronian triangles with diameter between $(x-1)\cdot 10\,000+1$ and $x\cdot 10\,000$.
The best least squares fitting function of the form $c_1+c_2\log x+c_3\log^2 x+c_4 x+c_5x\log x+c_6x\log^2 x$ is given by
$$
  160436.33
  +117761.45\log x
  +3191.78 \log^2 x
  +12023.76 x
  -2787.79 x\log x
  +169.14 x\log^2x 
$$
and leads to a $\Vert\cdot\Vert_2$-distance of $152\,331$ for the entire data.

%% The best least squares fitting function of the form $c_1+c_2x+c_3x^2+c_4\log x+c_5\log^2x+c_6\log^3x$ is given by
%% $$
%%  170065.682 - 86.3971x + 0.0323x^2 + 137440.26\log x - 1475.88\log^2 x + 2913.216\log^3 x
%%$$

%% As each million range of Heronian triangles were computed, a counting function was
%% determined, which fit the distribution.
%% 
%% The best least squares fitting function was found to be the sum of a quadratic polynomial in $x$ and a cubic polynomial in $\log(x)$:
%% \begin{equation}
%% count(x) = 170065.6818989618484 - 86.39708463701270545x + 0.03230092797640438240x^2 + 137440.2595431437350\log(x) - 1475.879585920842312\log(x)^2 + 2913.215936091200828\log(x)^3
%% \end{equation}
%% where $x$ is a 10K decade range for the Heronian triangle diameter, i.e. $x = 1$ is the range from 1 to 10,000, and $x = 443$ is the range from 4,420,001 to 4,430,000, etc.
%% 
%% \begin{center}
%%   \textbf{ToDo:} the equation needs to be split and a graph included
%% \end{center}

We remark that, besides the (implicit) $O(n^{1+\varepsilon})$ upper bound from \cite{herontriangles}, we
are not aware of any non-trivial lower and upper bounds for the number of (primitive) Heronian triangles
with a given diameter. As shown in \cite{1265.52018} one may deduce lower bounds for the minimum diameter of plane
integral point sets, where the current knowledge is still very weak \cite{1047.52011}, from such estimates.

\section{Exhaustive generation of $n$-clusters up to a given diameter}
\label{sec_exhaustive_n_cluster}

In order to determine the, with respect to its diameter, smallest $7$-cluster we have performed an
exhaustive search for $n$-clusters up to a given diameter. For the purpose of this paper the
chosen maximum diameter is $6\cdot 10^6$. A starting point is a complete list of all Heronian triangles
up to this diameter. More concretely we have chosen the exhaustive algorithm described in 
Section~\ref{sec_heronian_triangles} to generate all primitive Heronian triangles up to diameter
$6\cdot 10^6$ and extended this list by including all rescaled version such that the resulting diameter
is at most $6\cdot 10^6$.

The underlying basic idea to construct $n$-clusters is to combine two $n-1$-clusters sharing a common $n-2$
cluster. This way we can benefit from the fact that the constraints can be partially checked very early.
So starting from a list of $3$-clusters, i.e.\ Heronian triangles, we generate all $4$-clusters, then all
$5$-clusters, then all $6$-clusters, and finally all $7$-clusters.

For the first combination step, i.e.\ $n=4$, ``sharing a common $n-2$-cluster'' means that the two triangles
which should be combined both must have a side of the same length.

To avoid time-extensive duplicates and also the need to store extensive lists in memory we apply the concept
of orderly generation, see \cite{winner}, which avoid isomorphism search when cataloging combinatorial
configurations like in our example integral point sets or $n$-clusters. To this end a canonical form
has to be defined so that during the algorithm only canonical objects are combined. The constructed objects
are accepted if and only if they are canonical too. The benefit from such an approach is that no isomorphic
copies arise. For the details we refer the reader to \cite{1265.52018} with the adaptation
of considering triangles of characteristic $1$.

As a result we have computationally verified that the smallest $7$-cluster has diameter $2262000$ and that
there is no other $7$-cluster with diameter less then or equal to $4\cdot 10^6$. Along the way we have 
also exhaustively constructed all $4$-, $5$-, $6$-, and $7$-cluster with diameter at most $6\cdot 10^6$. Those
lists will be beneficial for the construction of additional $7$-clusters as will be explained in the 
following sections.

\section{Combining lists of $n$-clusters}
\label{sec_combine}
In the previous section we have described an algorithm to exhaustively generate a list of all $n$-clusters
up to given diameter taking a complete list of $n-1$-clusters with respect to that diameter. As induction
start we need a complete list of all Heronian triangles up to the used diameter. As described in 
Section~\ref{sec_heronian_triangles} the computational limits of such an approach force restrictions
to rather small diameters where only a few $7$-clusters exist. So from now on we will leave the approach
of exhaustive generation and go over to incomplete construction algorithms.

Our assumption for this section is that we are given a list of $n$-clusters, which we then combine to a list
of $n'$-clusters. For our paper, the most general setting is the following: Given a list $L_1$ of $n_1$-clusters and
a possibly different list $L_2$ of $n_2$-clusters we consider pairs $(l_1,l_2)$, where $l_1\in L_1$ and $l_2\in L_2$,
to construct $n'$-clusters, where mostly $n'>\max(n_1,n_2)$.

In Section~\ref{sec_exhaustive_n_cluster} we have assumed that the $n-1$-clusters $l_1$ and $l_2$ share a common
$n-2$-cluster. Since in the end we are only interested in lists of non-similar $n$-clusters we relax that to
the requirement that $l_1$ and $l_2$ contain a common $c$-cluster, where $c$ is an additional parameter.

Mostly we will restrict ourselves on the largest $c$-cluster of $l_1$ while looping over all non-isomorphic
$c$-clusters of $l_1$ is also possible. Having the $c$-cluster $C_1$ of $l_1$ fixed we loop over all $c$-clusters
$C_2$ of $l_2$ and check whether $C_1$ and $C_2$ can be rescaled so that they coincide. This check is implemented
as follows: Let $diam_1$ be the diameter of $C_1$ and $diam_2$ be the diameter of $C_2$ we define
$f_1=diam_2/gcd(diam_1,diam_2)$ and $f_2=diam_1/gcd(diam_1,diam_2)$. With this $C_1$ and $C_2$ are similar if and only
if $f_1\cdot C_1$ is isomorphic to $f_2\cdot C_2$. Comparing the sorted lists of the pairwise distances is a first
computationally cheap test for this task. If successful we compare the canonical forms of $C_1$ and $C_2$ and compare them.

So by rescaling we are in the situation that $l_1$ and $l_2$ contain a common $c$-cluster and we proceed by computing common
coordinates: We apply to the algorithm from Subsection~\ref{subsec_rational_coordinates} to compute
coordinates for $l_1$ and $l_2$ separately.\footnote{If $L_2$ is large it is computationally beneficial to store
a coordinate representation, given by the algorithm in Subsection~\ref{subsec_rational_coordinates}, for each $l_2\in L_2$. }
By assuming that the first $c$ points of $l_1$ and $l_2$ coincide we can
obtain a common coordinate system by just scaling the numerators. We remark that for $c=2$ we have two possibilities for the 
join, otherwise just one. Having the coordinates at hand we can look over all $k$-set of the points, where $k$ is sufficiently
large, and check whether they satisfy the conditions of a $k$-cluster, where we relax the condition of integral distances to
rational distances. If all (relaxed) conditions are satisfied we store a primitive version of the corresponding, possibly scaled,
$k$-cluster.

\bigskip

We have mostly used three instances of this general framework. The first is with the parameters $n_1=n-1$, $n_2=3$, and
$c=2$, i.e.\ we try to extend a given list of $n-1$-clusters by combining them with a list of primitive Heronian triangles
along a common edge. Since we use rescaling this combination is always possible, albeit it is not clear if there result
any $n$-clusters at all. Depending on the available computation time and the size of the list of the $n-1$-clusters one
may choose all known primitive Heronian triangles for the second list. We have done that to a large extend for the list
of known $7$-clusters but unfortunately did not locate an $8$-cluster. 

The second instance we used is with the parameters $n_1=n_2=6$ and $c=3$ to combine lists of $6$-clusters sharing a common triangle 
to obtain additional $6$- or $7$-clusters. The resulting point set consists of nine points. We remark that the second method was able
to discover some previously unknown $6$- and $7$-clusters but turned out to be rather slow. For later reference we call this method
the combine-hexagons algorithm.  Similar approaches seemed to be even less successful. 

The third method mimics the exhaustive generation method from Section~\ref{sec_exhaustive_n_cluster}, i.e.\
$n_1=n_2=n-1$ and $c=n-2$, starting from $n=4$ and increasing it by one in each iteration.

\section{Triangle extensions}
\label{sec_triangle_extension}
While the algorithms in Section~\ref{sec_combine} have to be iteratively used to end up with $n$-clusters for $n$ large, we will
now describe an algorithm that directly approaches $n$-clusters for $n$ as large as possible. Let $L$ be a list of primitive 
Heronian triangles of length $n$.

\begin{algorithm}
  {\textbf{(Triangle extensions)}\\}
  for $i$ from $1$ to $n$\\
  \hspace*{3mm}$\mathcal{P}=\emptyset$\\
  \hspace*{3mm}for $j$ from $i$ to $n$\\
  \hspace*{6mm}combine $L(i)$ with $L(j)$ in all possible ways\\ 
  \hspace*{9mm}compute coordinates of fourth point $p$\\
  \hspace*{9mm}if $L(i)\cup p$ is a $4$-cluster then add $p$ to $\mathcal{P}$\\
  \hspace*{3mm}compute all pairwise distances between the points in $\mathcal{P}$\\
  \hspace*{3mm}loop over all $k$-sets, where $k\ge 2$, $\mathcal{K}$  of $\mathcal{P}$\\
  \hspace*{6mm}if $L(i)\cup \mathcal{K}$ is a cluster and $k\ge 3$ then output $L(i)\cup \mathcal{K}$\\
\end{algorithm}

We remark that the loop over the $k$-sets is iteratively done, i.e.\ all subsets of the present $k$-set
are previously checked so that e.g.\ the upper bound for $k$ automatically is chosen. The implementation 
details for the coordinate and distances computations are described in Section~\ref{sec_low_level_details}. 

\section{Low level mathematical and implementation details}
\label{sec_low_level_details}
\noindent
In the previous sections we have described our algorithms without much
implementation details. Since the application of those algorithms
result in many sub computations like e.g.\ coordinate and distance
computations those sub routines have to be carefully designed in order
to save costly unlimited precision rational computations.

\subsection{Compute rational coordinates of a Heronian triangle}
\label{subsec_heron_coordinates}
Suppose we are given three integer side lengths $a$, $b$, and $c$, which
form a non-degenerate Heronian triangle. Our aim is to compute rational
coordinates for the points $P_1$, $P_2$, and $P_3$ attaining those pairwise
distances, i.e.\ $|P_1P_2|=a$, $|P_1P_3|=b$, and $|P_2P_3|=c$.

W.l.o.g.\ we can assume that the first point is located in the origin of 
our coordinate system, i.e.\ $P_1=(0,0)=\big(\frac{0}{2a},\frac{0}{2a}\big)$.
Considering the edge of length $a$, we locate the second point of the positive
side of the $x$-axis so that $P_2=(0,a)=\big(\frac{0}{2a},\frac{2a^2}{2a}\big)$.
By introducing variables $x$ and $y$ for the coordinates of the
third point, i.e.\ $P_3=(x,y)$, we obtain a quadratic equation system, which
can be uniquely solved if we assume that the $y$-coordinate is non-negative:
\begin{eqnarray*}
  x &=& \frac{b^2-c^2+a^2}{2a}=:\frac{t_1}{2a} \\
  y &=& \frac{\sqrt{4b^2a^2-(b^2-c^2+a^2)^2}}{2a}=\frac{\sqrt{4b^2a^2-t_1^2}}{2a}=:\frac{t_2}{2a}
\end{eqnarray*}
The second solution is given by $\big(\frac{t_1}{2a},-\frac{t_2}{2a}\big)$.

In some algorithms all permutations of the three edge lengths of a Heronian
triangle $(a,b,c)$ should be considered. To this end we assume that the above
auxiliary integer values $t_1$ and $t_2$ have already been computed. Permuting
the two latter side lengths, i.e.\ $(a,c,b)$, is equivalent to swap the points $P_1$
and $P_2$. The corresponding coordinates with non-negative $y$-values are given by
$$
  \Big(\frac{0}{2a},\frac{0}{2a}\Big),
  \Big(\frac{2a^2}{2a},\frac{0}{2a}\Big),
  \Big(\frac{2a^2-t_1}{2a},\frac{t_2}{2a}\Big)
$$
By applying a suitable rotation matrix we obtain the coordinate representation
$$
  \Big(\frac{0}{2b},\frac{0}{2b}\Big),
  \Big(\frac{2b^2}{2b},\frac{0}{2b}\Big),
  \Big(\frac{2b^2-t_1}{2b},\frac{t_2}{2b}\Big)
$$
for the triangle $(b,c,a)$ and
$$
  \Big(\frac{0}{2c},\frac{0}{2c}\Big),
  \Big(\frac{2c^2}{2c},\frac{0}{2c}\Big),
  \Big(\frac{2a^2-t_1}{2c},\frac{t_2}{2c}\Big)
$$
for the triangle $(c,a,b)$.

So there is no need to compute additional square roots. Of course the common subexpressions like
e.g.\ $a^2$, $b^2$, and $c^2$ should be stored additionally.

\subsection{Compute rational coordinates of an $n$-cluster}
\label{subsec_rational_coordinates}
We assume a suitable but fixed ordering of the $n$ points and denote the (integer) distance between the
first two points by $d$. W.l.o.g.\ the first point has coordinates $P_1=(0,0)=\big(\frac{0}{2d},\frac{0}{2d}\big)$
and the second $P_2=(0,d)=\big(\frac{0}{2d},\frac{2d^2}{2d}\big)$. For the third point we utilize the formula
for $(x,y)$ in Subsection~\ref{subsec_heron_coordinates}. We restrict ourselves onto the case of a positive $y$
coordinate. For points $4$ to $n$ we also compute the positive coordinates according to the previous 
Subsection~\ref{subsec_heron_coordinates} using the distances to point~1 and point~2. Next we have to check
whether the squared distance to the coordinates of point~3 coincides with the presumed squared distance and
possibly negate the computed $y$-coordinate. Thus all points have coordinates $\big(\frac{x_i}{2d},\frac{y_i}{2d}\big)$
with integers $x_i,y_i$.

\subsection{Checking for rational distances}
Suppose we are given two points with rational coordinates $\big(\frac{x_1}{a_1},\frac{y_1}{b_1}\big)$
and $\big(\frac{x_2}{a_2},\frac{y_2}{b_2}\big)$. The task is to decide whether they are at rational distance
and eventually compute the distance. Since during our searches most of the checked distances are irrational
it is important to have a quick check for the decision problem. An exact expression for the distance is given by
$$
  \frac{\sqrt{(b_1b_2)^2(a_2x_1-a_1x_2)^2+(a_1a_2)^2(b_2y_1+b_1y_2)^2}}{a_1a_2b_1b_2}.
$$
Thus the problem is reduced to the question whether a certain integer is a square.

Here we can benefit from modular arithmetic. Suppose that $m$ is an arbitrary integer and compute
$(b_1b_2)^2(a_2x_1-a_1x_2)^2+(a_1a_2)^2(b_2y_1+b_1y_2)^2 \mod m$ by performing all intermediate
computations modulo $m$. If the result is not a square in $\mathbb{Z}_m$ the distance under study
can not be rational. If $m$ is a product of distinct primes then we can check the square property
separately for each prime $p$ by simply tabulating a boolean incidence vector for the squares in 
$\mathbb{Z}_p$. In our implementation we use $m_1=493991355=3\cdot 5\cdot 11\cdot 13\cdot 17\cdot 19\cdot 23\cdot 31$
and $m_2=622368971=7\cdot 29\cdot 37\cdot 41\cdot 43\cdot 47$, i.e.\ we perform two successive modular
tests. Since computations modulo $4$ are very cheap in most arbitrary precision libraries it pays off
to first check whether the integer under study is equivalent to either $0$ or $1$ modulo $4$; otherwise
its square can not be rational.

If we can assume a common denominator of the coordinates, as e.g.\ implied by the algorithm in
Subsection~\ref{subsec_rational_coordinates}, the computations can be simplified since the distance
between the points $\big(\frac{x_1}{d},\frac{y_1}{d}\big)$ and $\big(\frac{x_2}{d},\frac{y_2}{d}\big)$ is given by
$$
  \frac{\sqrt{(x_1-x_2)^2+(y_1+y_2)^2}}{d}.
$$
\subsection{Canonical forms}
In order to be able to check $n$-clusters for similarity we define a canonical form
in such a way that two $n$-clusters are similar if and only if their canonical forms
coincide. Given a matrix of the pairwise rational distances we first normalize by
multiplying with the least common multiple of the denominators and then by 
dividing the greatest common divisor of the resulting nominators. Now we are given 
integer distances whose greatest common divisor is trivial. Since distances are symmetric
it suffices to consider the upper right triangular submatrix without the diagonal of zeros.
Appending the columns of this matrix gives a vector -- distance vector for brevity, which
can be compared with respect to the lexicographical ordering. We define  the canonical
form to be the lexicographically maximal distance vector over all permutations of the points.

Clearly the defined canonical form is unique and we can determine it by comparing all $n!$
possible permutations. For our purposes this was fast enough even for $n=7$, but we remark 
that one can easily design $O(n^3)$ algorithms to compute the canonical form.

\section{Circle inversion}
\label{sec_circle_inversion}
As observed in \cite{1209.52009} the rationality of distances in $\mathbb{R}^2$ is preserved by translations,
rotations, scaling with rational numbers and by some kind of circle inversion. We go into the details of the
latter transform. Assume that our point set has a point in the origin, then a circle inversion through the origin
with radius one sends each point with coordinates $(x,y)$ besides the origin to
$\big(\frac{x}{x^2+y^2},\frac{y}{x^2+y^2}\big)$\footnote{Using complex notation this is (ignoring a reflection)
equivalent to the map $z\mapsto\frac{1}{z}$.}.

Using this transform we can construct $n-1$-clusters from $n$-clusters by moving each of their points to the
origin and applying the described circle inversion. Doing this for the set of all known $7$-clusters gives
no new $6$-clusters, while even preserving the set of the contained subtriangles, i.e.\ the set of the (normalized)
subtriangles from the resulting $6$-clusters coincides with the set of the subtriangles contained in the $7$-clusters.

\highlightt{Discarding one point is, on the one hand disadvantageous, but gives us some freedom in the initial point 
set, it does not have to be an $n$-cluster. To be more precise, we need rational point sets $\mathcal{P}$ with
characteristic~$1$, where no three points are on a line and no four points are on a circle.} 
%% Discard one point is on the one hand disadvantageous but gives us some freedom in the initial point set, i.e.\ it
%% does not have to be an $n$-cluster. To be more precisely we need rational point sets $\mathcal{P}$ with characteristic
%% one, where no four points are on a line and no four points are on a circle. 
Circle inversion at a vertex of $\mathcal{P}$ automatically
destroys the property of those collinear triples. We were able to extend some of the $7$-clusters to $8$-point rational sets.
Clearly at most $3$ points are on a line and all those lines intersect in the $8$th point. Unfortunately in each of this cases
the $8$th point also was part of a circle containing four points of the point set. A promising configuration might be the 
so-called Pappus configuration consisting of nine points and nine line, with three points per line and three lines through
each point. Unfortunately we were not able to find a representation of the Pappus configuration with pairwise rational distances.

So while circle inversion might be theoretically interesting we were not able to draw any computational advantages.

\section{Points at rational distance from the vertices of a triangle}
\label{sec_fourth_point}

\highlight{
Instead of extending a given $(n-1)$-cluster with the aid of Heronian triangles one might directly appeal
to Almerings theorem that the set of points at rational distance to the vertices of a sub triangle of a cluster
is dense in the plane. Here we simplify notation and assume that we are given a triangle with rational side lengths.
Moreover we assume that this triangle is rotated into a convenient position, compare Subsection~\ref{subsec_heron_coordinates}.
In the following two subsections we present two approaches capable of producing candidates for a forth point
with rational distances to the vertices of the initial triangle. We have applied both methods in order to
extend $n$-clusters for small $n$ and report that they both find some examples but generally the necessary running
time, i.e.\, the number of choices for the method's parameters, is not competitive compared to the other algorithms
described earlier on. 
}

\highlight{
\subsection{Pythagorean arctangent method}
Given an arbitrary Heronian triangle with side lengths $a$, $b$, and $c$ we rotate it in the convenient
position shown in Figure~\ref{fig:Heron_tangents}, i.e., the three vertices have coordinates $A=(0,0)$, 
$B=(c,0)$, and $C=\left(\frac{b^2+c^2-a^2}{2c},\frac{2A_\Delta(a,b,c)}{c}\right)$, where 
$A_\Delta(a,b,c)\in\mathbb{N}\subset\mathbb{Q}$ is the area of the triangle. The method of {\it Pythagorean arctangent}
uses two Pythagorean arctangents to choose two smaller angles, one $\theta$ located at $A$
and the other $\phi$ at $B$ such that the intersection point $D$ has rational coordinates
and the distances $\overline{AD}$ and $\overline{BD}$ are rational too.  In other words the triangle 
spanned by vertices $A$, $D$, and $B$ is made Heronian.
\begin{figure}[htp]
 \begin{center}
 \includegraphics[width=4.0in]{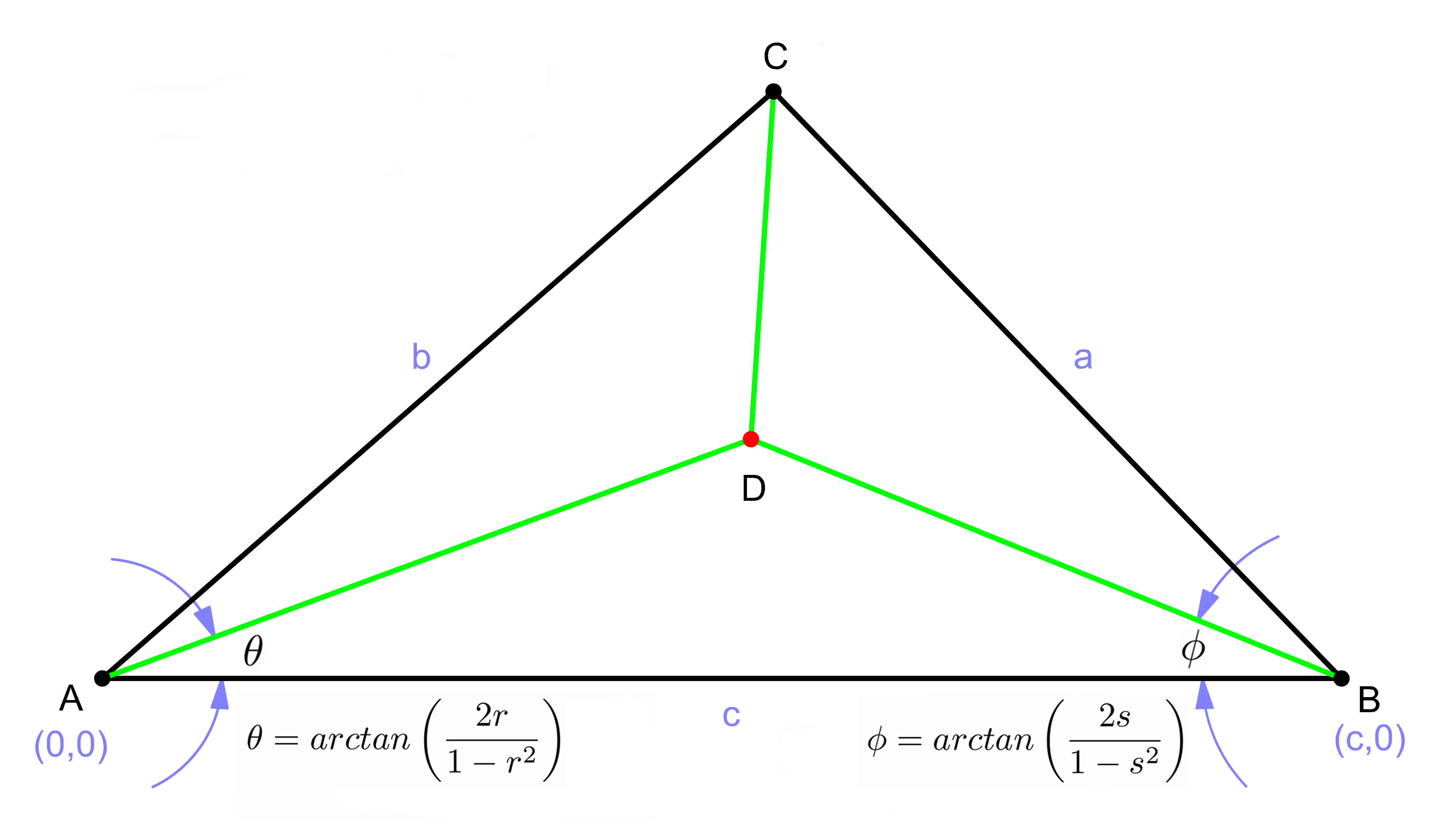}
 \caption{Pythagorean arctangent method}
 \label{fig:Heron_tangents}
 \end{center} 
\end{figure}
A Pythagorean arctangent can be parameterized as $\frac{2mn}{(m+n)(m-n)}$ with integers $m$ and $n$.
Dividing both the numerator and the denominator by $m^2$ we obtain the rational parameterizations
$\theta=\arctan\!\left(\frac{2r}{1-r^2}\right)$ and $\phi=\arctan\!\left(\frac{2s}{1-s^2}\right)$ with $r,s\in(0,1)$.
With this we can easily check that the coordinate of the fourth point
$$
  D=\left(\dfrac{c s  \left(1-r^2\right)}{r + s - rs\left(r+s\right) } \;,\; \dfrac{2 c r s}{r + s - rs\left(r+s\right) }\right)
$$
and also the two distances
$$
  \overline{AD} = \dfrac{ c  s  \left(1+r^2\right) } { r s\left(r+s\right) -r -s }
$$
and
$$
  \overline{BD} = \dfrac{ c  r  \left(1+s^2\right) } { r s\left(r+s\right) -r -s }
$$
are rational. For the third distance we obtain
$$
  \overline{CD} = \sqrt{ \left( \dfrac{b^2+c^2-a^2}{2c} - \dfrac{cs\left(1-r^2\right)}{r+s-rs\left(r+s\right)} \right)^2 
  \!\!\!+\!\! \left(\dfrac{A_\Delta(a,b,c)}{\frac{1}{2}\cdot c} - \dfrac{2crs}{r+s-rs\left(r+s\right)} \right)^2 }.
$$
Using the substitution $Y = r + s - rs\left(r+s\right)$ yields the simplified expression
$$
  \overline{CD} = \dfrac{ \sqrt{ \left( Y\left(b^2+c^2-a^2\right)-2c^{2}s\left(1-r^2\right)\right)^2 + \left( 4A_\Delta(a,b,c) Y - 4c^{2}rs \right)^2 } } { 2cY}
$$
The numerator of this equation must be a square, if $\overline{CD}$ has to be rational.  Unfortunately it has the form of a
homogeneous quartic equation in $r$, $s$ which, in general, is difficult to solve.  If one point for a choice of $r$, $s$ is
found, then others exist, because it can be transformed into an elliptic curve. At this time, the only choice is to actually
determine the value and then take the square root for arbitrary choices of $r$, $s$. $Y$ is fixed by our choice of $r$, $s$,
but the numerator-constraint must be taken into account.
} 

\highlight{
\subsection{Exploiting Ceva's theorem}
Given a triangle with vertices $A$, $B$, and $C$, let the lines $AO$, $BO$ and $CO$ be drawn,
where $O$ is a common point. The intersection points at the sides of the initial triangles are denoted
by $D$, $E$, and $F$, respectively, see Figure~\ref{fig:Ceva_point}, where $O$ is the common intersection
point. Let $|XY|$ denote the signed length of the segment between $X$ and $Y$, i.e., 
$|YX|=-|XY|$. With this notation, Ceva's theorem states
$$\frac{|AF|}{|FB|}  \cdot \frac{|BD|}{|DC|} \cdot \frac{|CE|}{|EA|} = 1.$$
\begin{figure}[htp]
 \begin{center}
 \includegraphics[width=3.0in]{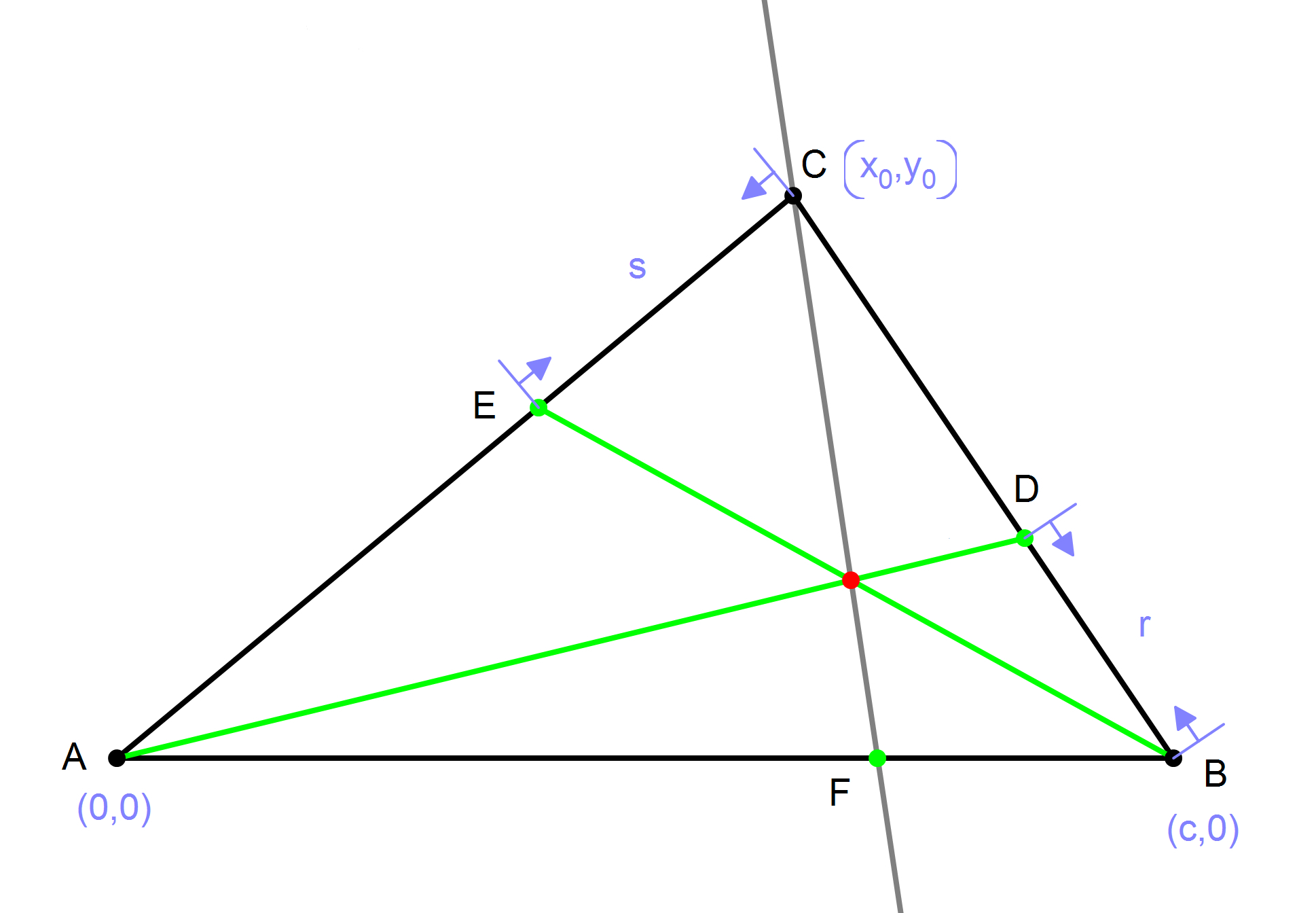}
 \caption{General Ceva Point}
 \label{fig:Ceva_point}
 \end{center}
\end{figure}
We set $r=\overline{DB}/\overline{BC}$, $s=\overline{EC}/\overline{AC}$ and use the rational coordinates $x_0,y_0$
of $C$ to state the coordinates
$$
  \left (\dfrac{(1-s)\cdot \left (-c+c\cdot r-r\cdot x_{0}\right )}{-1+s-r\cdot s},\dfrac{-r\cdot y_{0}\cdot (1-s)}{-1+s-r\cdot s}\right )
$$
of the common intersection point $O$, which clearly are rational. For the three distances we obtain
\begin{eqnarray*}
  \overline{AO} &=& \dfrac{ \sqrt{ \left( c\left(1-r\right) + r\cdot x_{0} \right) ^2 + \left(r\cdot y_{0}\right)^2 }  \cdot \left| s-1 \right| }{\left| -1 + s - r \cdot s \: \right| },\\
  \overline{BO} &=& \dfrac{\sqrt{ \left( c - x_{0}\cdot\left(1-s\right)\right)^2 + \left( y_{0}\cdot \left( 1-s \right) \right)^2 } \cdot \left| r \right| }{\left| -1 + s - r \cdot s \: \right| },\text{ and}\\
  \overline{CO} &=& \sqrt{ \left( x_{0} - \dfrac{\left(r\cdot x_{0} + c\left(1-r\right)\right) \cdot \left(s-1\right)}{-1+s-r\cdot s} \right)^2 + \left( y_{0} - \dfrac{r\cdot y_{0} \cdot\left(s-1\right)}{-1+s-r\cdot s} \right)^2 }.
\end{eqnarray*}
We can easily constrain $r$ and $s$ so that the first two distances get rational. Having chosen suitable rational numbers
$r,s\in(0,1)$, we can then check the third distance $\overline{CO}$. Viewed from a different angle we choose
a Heronian triangle with side lengths $a$, $b$, $c$ and another one with side lengths $d$, $e$, $f$. Via scaling we ensure
that $a$ and $d$ coincide. It remains to check that the two off-axis points ($C$ and $D$ in Figure~\ref{fig:Ceva_point}) 
are at a rational distance -- compare Section~\ref{sec_triangle_extension}.
}

\section{Choosing promising Heronian triangles}
\label{sec_promising}
The algorithms presented in the previous sections can in principle deal with large lists of $n$-clusters, but
of course the computation time limits such searches. In order to find many non-similar $7$-clusters we have 
tried to restrict ourselves on promising search spaces. Either the exhaustive-like algorithm from Section~\ref{sec_combine}
or the triangle extension algorithm from Section~\ref{sec_triangle_extension} grounds on a list of Heronian
triangles and then dives into the resulting search space. Unfortunately we do not have the computational capacity
to start those algorithms with all Heronian triangles known to us but have to select a subset of them. Of course
this subset should be selected in a way so that is small on the one hand but generates many $7$-clusters. To
satisfy the latter aim is essential but of course the hardest part. So far we have no theoretical justification 
but only computational results in that direction.  Conceptually the best way is to invent a method that is able
to compute a score for a given Heronian triangles and then choosing a given number of Heronian triangles with
the largest score. 

A very easy but effective scoring function is the negative diameter of all  Heronian triangles. In order to verify our
claim we have used the triangle extension algorithm with subsets of $1000$ Heronian triangles. Using the first $1000$
smallest, with respect to diameter, Heronian triangles produces $237$ $6$-clusters and four $7$-clusters (having diameters
$5\,348\,064$, $15\,772\,770$, $47\,570\,250$, and   $662\,026\,750$). The second smallest $1000$ Heronian triangles 
produces \textit{only} nine $6$-clusters and no $7$-cluster. 

\medskip

A promising idea might be to use the number of divisors or prime divisors of the side lengths normalized by
magnitude, i.e.\ prime side lengths should get the lowest possible score while highly composite numbers get 
large scores. Exemplary we report the results of two explicit scoring functions based on this idea.
For $$score_1(a,b,c):=\frac{\#\text{prime divisors }a}{\log \log a}+\frac{\#\text{prime divisors }b}{\log \log b}+
\frac{\#\text{prime divisors }c}{\log \log c}$$ we have chosen the $1000$ Heronian triangles with maximal score
among all Heronian triangles with diameter at most $10\,000$. Applying the triangle-extension algorithm results in 
three $6$-clusters and no $7$-cluster.
The similar function $$score_2(a,b,c):=\frac{\#\text{prime divisors }a}{\log a}+\frac{\#\text{prime divisors }b}{\log b}+
\frac{\#\text{prime divisors }c}{\log c}$$
increases the number of found $6$-clusters to $40$ within the same setting. But of course $score_2$ tends to prefer
triangles with smaller diameter. We remark that using the number of divisors instead of the number of prime divisors
yields similar results.

\medskip

The most successful approach in our computationally study was to use the known lists of $n$-clusters as selectors.
To be more precisely, given a list of $n$-clusters we can determine the contained sub-triangles, which then,
after rescaling, gives a list of primitive Heronian triangles. If the resulting list of Heronian triangles is too
large for our purposes we take the $m$ smallest ones according to their diameter or we take frequency into account,
i.e.\ we consider only those primitive Heronian triangles which appear at least $k$ times, where $k$ is
suitably chosen, as sub-triangles within the list of $n$-clusters.

Exemplary we report the following to experiments performed near the end of our computational study, were we already know
lots of $6$- and $7$-clusters. For $n=6$ and $n=7$ we choose the $1000$ Heronian triangles having the smallest diameter, 
respectively. In the first case triangle extension yields $247$ $6$-clusters and four $7$-clusters. For the latter case
we obtain $912$ $6$-clusters and $100$ $7$-clusters. So a higher initial value of $n$ results in more clusters, but of course
those examples are harder to find.

\medskip

A completely different idea is to associate Heronian triangles $(a,b,c)$ with ellipses represented by $\frac{a+b}{c}$.
As an experiment we took the $3\,000\,000$ smallest Heronian triangles and computed the three associated ellipses in each case.
The most frequent ellipse representation occurs $10\,277$ times. Taking the smallest $1000$ triangles results in $603$ $5$-clusters
applying the triangles extension algorithm. Taking triangles from ellipse representations that occur exactly once result in just 
six $5$-clusters.

\section{Computational observations}
\label{sec_observations}

In this section we collect some computational observations that help us
to design our searches for $7$-clusters.

\begin{claim}
  The triangle-extension algorithm is more effective than the combine-hexagons algorithm. 
\end{claim}

Using the $412$ triangles contained in the original twenty-five $7$-clusters found by Simmons and Noll in 2010
as an input for the triangle-extension algorithm yields $84$ non-similar $7$-clusters in less than two minutes
computation time. If we instead take the sub-hexagons of the original twenty-five $7$-clusters plus an additional
list of $1736$ hexagons and apply the combine-hexagons algorithm we end up in $33$ non-similar $7$-clusters. We remark
that all but one of these heptagons is contained in the list of the $84$ heptagons from the triangle extension algorithm.
Additionally the computation time of the combine-hexagons algorithm is usually much larger than the computation time
of the triangle-extension algorithm.

\begin{claim}
  Stripping isosceles triangles from the input set of Heronian triangles
  only mildly reduces the number of $6$- and $7$-clusters found in the search
  of the triangle-extension algorithm.
\end{claim}

Because any pair of isosceles Heronian triangles forms a $4$-cluster, there
are numerous $4$-clusters formed from pairs of isosceles triangles.  When
testing a pair of $4$-clusters that have this property, the pair will not be
interesting because three points will lie on the line through the median of
the base of the isosceles triangles. 

As expected the runtime increases while including isosceles Heronian triangles, where
the precise factor heavily depends on the chosen subset of Heronian triangles. For comparison we
have chosen the $1000$ smallest non-isosceles Heronian triangles and applied the triangle-extension
algorithm, which resulted in $172$ $6$-clusters and four $7$-clusters.  So we have missed $65$
$6$-clusters but no $7$-cluster. Here the computation time was decreased by a factor of two.
In a larger experiment we have chosen $1\,383\,799$ Heronian triangles and obtained $424\,593$ $6$-clusters
and $1\,110$ $7$-clusters. Stripping all $24\,583$ isosceles triangles we have obtained
$424\,543$ $6$-clusters and $1\,110$ $7$-clusters, while the computation time decreases by a factor
larger than $10$.

\begin{claim}
  \label{claim_partionioning}
  Partitioning the set of triangles can speed up the search of the triangle-extension algorithm.
\end{claim}

Given a list of $m$ $n$-clusters containing the same $n-1$-cluster the ordinary combination would need
$m^2$ tests. Since integral point sets with many points on a line or a circle are quite common
it makes sense to take this fact into account. Partitioning $4$-clusters by a line through $2$ of
the points or by a circle through $3$ of the points avoids many spurious comparisons and speeds up the search.
The important thing is that a pair of items in a partition cannot form an $n+1$-cluster because it would
violate a con-circularity or co-linearity constraint. In our programs we can either turn on and off the
partitioning algorithm, but mostly use it to increase the computation speed. The typical performance boost 
is around 10~\%.% for medium- and large-sized experiments.

\begin{claim}
  \label{claim_no_large_triangles}
  Large Heronian triangles tend to not form 4-clusters.
\end{claim}

That is, given two random small Heronian triangles, the probability they form a $4$-cluster is
relatively high compared to the probability that two large Heronian triangles will form a $4$-cluster, i.e.\
we have to perform many unsuccessful combinations of Heronian triangles per found $4$-cluster.
To additionally justify this theoretically, one might appeal to Ceva's theorem. As we allow the size of a Heronian triangle
to increase the prime factors present in the numerators of the sines of the Heronian angles increase making
it more difficult to find sets of angles where the numerators cancel each other out. 

\begin{claim}
  Iterating the triangle-extension algorithm can find new triangles and $n$-clusters.
\end{claim}

As described in Section~\ref{sec_promising} combining the triangles contained in the twenty-five
$7$-clusters found by Simmons and Noll in 2010 yields $84$~non-similar $7$-clusters. Those $7$-clusters
contain $602$ triangles which combine to $86$ non-similar $7$-clusters using the triangle extension algorithm. 
Then the iteration gets stuck since those $7$-clusters contain exactly $602$ non-similar triangles again.

Similarly we have used the $237$~$6$-clusters which arose from combining the $1000$ smallest Heronian triangles,
see Section~\ref{sec_promising}. Those $6$-clusters contain $1808$ non-similar triangles which can be combined to 
$1644$ non-similar $6$-clusters and $22$ non-similar $7$-clusters.

\begin{claim}
  The rational distance test rules out most of the combinations of Heronian triangles.
\end{claim}

To verify this claim we report the statistics of a large scale experiment. We have chosen the $3\,000\,000$ smallest
primitive Heronian triangles along with those contained in the $6$-clusters known to us. Using $25\,000$ cores during
$4.5$ days $3.0~\cdot 10^{14}$ pairs of $3$-clusters were tried. In $99.71$~\% the missing sixth distance was not
rational. The concircular test ruled out $10\,414\,450\,261$ possibilities ($0.00$~\%) and the collinearity test
$20\,129\,596\,307$ possibilities ($0.01$~\%), while we found $835\,620\,202\,676$ (possibly similar) successful combinations
($0.28$~\%). The longest list of $4$-clusters containing a common $3$-cluster had length $396\,442$.
In Table~\ref{table_checks} we have summarized the corresponding statistics for the combinations of the resulting $k$-clusters
for $3\le k\le 7$.

\begin{table}[htp]
\begin{center}
\begin{tabular}{rrrrrrr}
  \hline
  $k$ & comb.\ & distance & concircularity & collinearity & successful & intersectable \\
  3 & $3.0~\cdot 10^{14}$ & $99.71$~\% & $0.00$~\% & $0.01$~\% & $0.28$~\% & 396442 \\
  4 & $2.1~\cdot 10^{15}$ & $41.87$~\% & $58.13$~\% & $0.00$~\% & $0.00$~\% & 91 \\
  5 & $1.6\cdot 10^8$ & $49.93$~\% & $33.17$~\% & $14.01$~\% & $2.89$~\% & 16 \\
  6 & $1.5\cdot 10^5$ & $60.89$~\% & $18.93$~\% & $8.86$~\% & $11.32$~\% & 2 \\
  7 & $82$ & $100$~\% & $0.00$~\% & $0.00$~\% & $0.00$~\% & 0 \\
  \hline
\end{tabular}
\caption{Failure of different checks for $k+1$-clusters combining two $k$-clusters}
\label{table_checks}
\end{center}
\end{table}

\section{Computational results}
\label{sec_results}

We have constructed \numsevenclusters~non-similar $7$-clusters and  \numsixclusters~non-similar $6$-clusters\footnote{The list
of the primitive $6$- and $7$-clusters currently known to us can be obtained at \cite{HP}.}.
The $5$- and $4$-clusters are so numerous that we did not collect them. The total number of stored Heronian triangles is $807\,677\,361$.
The smallest diameter of a primitive $7$-cluster is $2262000$ while the largest found primitive $7$-cluster has a diameter of 
$$92986018038515228913684944937313015456\approx 10^{38}.$$ The \numsevenclusters~$7$-clusters contain in total 
${7 \choose 3}\cdot 1154=40390$ sub-triangles, while only $9264$ of them are non-similar, i.e., on average each
(normalized) triangle is used more than four times. The smallest contained triangle is $(5,4,3)$, which is indeed the smallest
possible Heronian triangle, and the largest has diameter $121990813408205791\approx 10^{18}$. Some counts of $7$-clusters are
given in Table~\ref{table_seven_cluster}. We remark that the Heronian triangles $(6,5,5)$, $(8,5,5)$, and $(13,12,5)$ are not
contained in any of the known $7$-clusters. The $6$-clusters contain more than $1\,400\,000$ non-similar Heronian triangles.
The smallest Heronian triangle that is not contained in one of the known $6$-clusters is $(149,148,3)$.

\begin{table}[htp]
  \begin{center}
    \begin{tabular}{rrrrrr}
      \hline
        diameter & \# $7$-clusters &  diameter & \# $7$-clusters &  diameter & \# $7$-clusters \\
        $\le 10^7$    &   4 & $\le 10^{19}$ &  688 & $\le 10^{31}$ & 1130 \\
        $\le 10^8$    &  11 & $\le 10^{20}$ &  752 & $\le 10^{32}$ & 1137 \\
        $\le 10^9$    &  26 & $\le 10^{21}$ &  819 & $\le 10^{33}$ & 1145 \\
        $\le 10^{10}$ &  52 & $\le 10^{22}$ &  877 & $\le 10^{34}$ & 1147 \\
        $\le 10^{11}$ &  89 & $\le 10^{23}$ &  927 & $\le 10^{35}$ & 1150 \\
        $\le 10^{12}$ & 139 & $\le 10^{24}$ &  974 & $\le 10^{36}$ & 1153 \\
        $\le 10^{13}$ & 198 & $\le 10^{25}$ & 1024 & $\le 10^{37}$ & 1153 \\
        $\le 10^{14}$ & 270 & $\le 10^{26}$ & 1050 & $\le 10^{38}$ & 1154 \\
        $\le 10^{15}$ & 347 & $\le 10^{27}$ & 1067 \\
        $\le 10^{16}$ & 431 & $\le 10^{28}$ & 1087 \\
        $\le 10^{17}$ & 516 & $\le 10^{29}$ & 1111 \\
        $\le 10^{18}$ & 609 & $\le 10^{30}$ & 1124 \\
      \hline
    \end{tabular}
    \caption{Number of (known) non-similar $7$-clusters up to a given diameter}
    \label{table_seven_cluster}
  \end{center}
\end{table}

As hardware we have used $25\,000$ cores at Google Inc.\ and the Linux computing cluster of the University of Bayreuth,
which consists of 201~2xIntel E5520 2.26 GHz and  52~2xIntel E5620 2.4GHz processors (100-300 jobs are done in parallel).
The computations for the triangle-extension algorithm using the triangles in the known $7$-clusters were done on a customary
laptop computer in less than one day of computation time per iteration. We used the GNU MP Bignum
library\footnote{http://gmplib.org/} and class library of numbers (CLN)\footnote{http://www.ginac.de/CLN/} libraries to
provide arbitrary precision integers and rationals. 

Although we have invested a large amount of processing power over the past year in our
current collaboration, we have not found an $8$-cluster.

\section{Conclusion}
\label{sec_conclusion}

\noindent
The techniques of finding $n$-clusters have dramatically improved since the
discovery of the first $6$-clusters in $\mathbb{R}^2$: back when some researchers
incorrectly conjectured that $6$-clusters in $\mathbb{R}^2$ did not exist.
At the current state it is still a significant computational challenge to
find new $7$-clusters, but we have shown that quite some examples exist.
A toolbox of algorithms to generate $n$-clusters is provided. Using the
triangle-extension algorithm one may eventually extend a small list of
$n$-clusters to a larger list of $n$-clusters by just combining their
contained sub triangles. Compared with its running time and its output
in terms of newly found $n$-clusters this is certainly the most effective
algorithm that is currently known. For a given $n$-cluster the knowledge
of only $n-2$ of its sub triangles may suffice to recover all distances and
so all ${n \choose 3}$ sub triangles. Moreover we have some kind of scale
invariance, i.e.\ only the angles but not the side lengths have to be known
in advance. Considering all possible scalings comes at constant cost.

However this algorithm is indentured from a good list of Heronian triangles,
or indirectly a list of starting $n$-clusters. To some extent the algorithm
itself produces some new Heronian triangles so that it can be applied iteratively.
But admittedly the number of successful iterations is observed to be rather small
in practice. So different algorithms are needed to populate the set of
\textit{promising} triangles. Choosing them directly from the list of Heronian
triangles, based on a scoring function, still has no satisfactory  solution
and is left as an open problem. So still the discovery of new $7$-clusters
depends on extensive computer calculations to that highly optimized
low level routines are essential to check a large number of cases.   

Along the way we have exhaustively constructed all primitive Heronian triangles
with diameter up to $6\cdot 10^6$. This database may serve as a starting point
to check various conjectures.

The question whether there exists an infinite number of non-similar $7$-clusters
is still open. At this point we would be the last to speculate that there are
anything but an infinite number of $8$-clusters in $\mathbb{R}^2$. 

%\bibliographystyle{amsplain}
%\bibliography{7-cluster}

\end{document}